\documentclass{amsart}

\DeclareSymbolFont{AMSb}{U}{msb}{m}{n}
\DeclareMathSymbol{\subsetneq}{\mathrel}{AMSb}{"28}

\newtheorem{theorem}{Theorem}[section]

\newtheorem{definition}[theorem]{Definition}
\newtheorem{lemma}[theorem]{Lemma}

\newtheorem{remark}[theorem]{Remark}

\newcommand{\N}{\mathbb{N}}
\newcommand{\Z}{\mathbb{Z}}

\newcommand{\AW}{\widetilde{W}}
\newcommand{\ZZ}{Z}
\newcommand{\Span}{\text{span\,}}
\newcommand{\sgn}{\text{sgn}}
\newcommand{\infdi}{I_2(\infty)}

\begin{document}
\title{The finite antichain property in Coxeter groups}
\author{Axel Hultman}
\address{Department of Mathematics, KTH, SE-100 44 Stockholm, Sweden }
\email{axel@math.kth.se}

\begin{abstract}
We prove that the weak order on an infinite Coxeter group contains
infinite antichains if and only if the group is not affine.
\end{abstract}

\maketitle

\section{Introduction}
Given an infinite poset, a natural problem is to decide whether or
not it contains infinite antichains (sets of pairwise
incomparable elements). It is known \cite{bjorner,higman} that every
antichain in the Bruhat order on any Coxeter group is
finite. Here, we consider the other of the two most common ways to
order a Coxeter group, namely the weak order. One observes that the
answer must depend on the group; it is straightforward to check that
the antichains in the infinite dihedral group are finite,
whereas there are infinite ones in the universal Coxeter group of rank $3$. The
open problem of characterising the groups with infinite antichains is
\cite[Exercise 3.11]{BB}. The main result of this paper is that the answer is
the following new characterisation of affine Weyl groups:
\begin{theorem}\label{th:main}
The weak order on an infinite Coxeter group contains an infinite
antichain if and only if the group is not affine.
\end{theorem}

After establishing notation in Section \ref{se:prel}, we use the
remaining two sections to prove our main result. In Section
\ref{se:general}, we show that affine groups do not possess infinite
antichains, whereas (irreducible) not locally finite ones do. The groups that remain
are the compact hyperbolic Coxeter groups. In Section \ref{se:comphyp}, it is
shown that they all have infinite antichains, thereby finishing the proof of
Theorem \ref{th:main}. While the proofs in Section \ref{se:general}
are uniform, we have been forced to resort to a case-by-case argument
in Section \ref{se:comphyp}.

\begin{remark}
{\em The Poincar\'e series of an affine
  Weyl group is given by a simple formula first proved by Bott
  \cite{bott}, see \cite[Theorem 8.9]{humphreys}. From it, it
  immediately follows that every affine Weyl group except the infinite
  dihedral group $\infdi$ has the following property: the number of
  elements of Coxeter length $k$ grows 
arbitrarily large as $k$ tends to infinity. Distinct elements of
the same length are always incomparable under the weak
order. Therefore, the only infinite Coxeter groups with bounded antichain
  size are $W\times \infdi$, for finite $W$. All other infinite groups
  have arbitrarily large finite antichains.}
\end{remark}

\section{Preliminaries}\label{se:prel}
We assume the reader to be familiar with basic theory of Coxeter groups
and root systems as can be found e.g.\ in \cite{BB} or
\cite{humphreys}. Here, we review scattered pieces of the
theory in order to agree on notation. For the most part, we borrow our
terminology from \cite{BB}. 

Throughout the paper, $(W,S)$ will denote a Coxeter system with
$|S|<\infty$. Given $J\subseteq S$, $W_J = \langle J\rangle$ is the
parabolic subgroup generated by $J$. Every coset in $W/W_J$ has a unique
representative of minimal length; the set of such representatives is
denoted by $W^J$.

We use $\ell(w)$ to denote the Coxeter length of $w\in W$. The
{\em right descent set} of $w$ is 
\[
D_R(w) = \{s\in S\mid \ell(ws) < \ell(w)\}. 
\]

\begin{definition}
The {\em (right) weak order} on $W$ is defined by $v\leq_Rw$ if and only
if there exists $u\in W$ such that $w = vu$ and $\ell(w) = \ell(v)+\ell(u)$.
\end{definition}
One can also define the left weak order $\leq_L$ in the obvious
way. In this paper, ``weak order'' always refers to the right weak
order. The results are of course equally valid for the left version.

Every group element $w\in W$ can be expressed as a word in the free
monoid $S^*$. If such a word has length $\ell(w)$ it is called a {\em
  reduced expression}. The combinatorics of reduced expressions is a
key to many properties of Coxeter groups and plays a prominent role in
our arguments. Abusing notation, we will sometimes blur the
distinction between elements of $W$ and their representatives in
$S^*$. We trust the context to make the meaning clear.

For $s,s^\prime\in S$, let $m(s,s^\prime)$ denote the order of
$ss^\prime$. This information is collected in the
{\em Coxeter diagram} which is a complete graph on the vertex set $S$
in which the edge $\{s,s^\prime\}$ is labelled with
$m(s,s^\prime)$. For convenience, we agree to suppress
edges with label $2$ and labels that are equal to $3$. 

Consider a word in $S^*$ representing $w\in W$. Deleting a factor
$ss$, we obtain another word representing 
$w$. Similarly, replacing a factor $ss^\prime ss^\prime s\dots$ of
length $m(s,s^\prime)$ with the factor $s^\prime ss^\prime ss^\prime
\dots$ of the same length, we again obtain a representative of
$w$. The former operation is called a {\em nil move}, the latter a
     {\em braid move}. Obviously, nil moves can never be performed on
     reduced expressions. The following important result is due to
     Tits \cite{tits}.  
\begin{theorem}[Word Property]$~$
\begin{itemize}
\item[{\rm (a)}]
Any word in $S^*$ can be brought to a reduced expression by a
sequence of braid moves and nil moves.
\item[{\rm (b)}] Given a reduced expression for $w$, every other reduced
  expression for $w$ can be obtained by a sequence of braid moves.
\end{itemize}
\end{theorem}
One consequence of the Word Property is that every reduced expression
for $w\in W$ uses the same set of generators. We use
$S(w)\subseteq S$ to denote this set.

\subsection{A recognising automaton for reduced expressions}
It is a fundamental fact that the language of reduced expressions in
$W$ is regular, i.e.\ recognised by a finite state automaton. In other
words, there is a directed graph on a finite vertex set whose edges
are labelled with elements from $S$, such that the sequences of labels
along directed paths beginning in some distinguished starting vertex
are exactly the reduced 
expressions for elements in $W$. The existence of such an automaton is
essentially due to Brink and Howlett \cite{BH}. They in fact showed
that the language of so-called normal forms is regular, but this is
enough due to a result of Davis and Shapiro \cite{DS}.

At one point in Section \ref{se:comphyp}, we will rely on explicit
computations in a particular recognising automaton for reduced
expressions. For this purpose, we briefly sketch how the automaton
works. For the sake of brevity, the reader will be kept on a
need-to-know basis. A more thorough account of the construction can be
found in \cite[Section 4.8]{BB}, which is largely based on material from
the thesis of Eriksson \cite{eriksson}. See also Headley's thesis
\cite{headley}. 

Suppose $\Phi$ is a root system for $W$ with simple roots
$\Delta = \{\alpha_s\mid s\in S\}$. A symmetric bilinear form on $V =
\Span \Phi = \Span \Delta$ is defined by
\[
(\alpha_s\mid\alpha_{s^\prime}) = -\cos \frac{\pi}{m(s,s^\prime)}.
\]
Given $w\in W$, we
recursively define a corresponding set $D_\Sigma(w)\subseteq \Phi^+$
of positive roots by
\[
D_\Sigma(e) = \emptyset
\]
(where $e\in W$ is the identity element) and, if $s\not \in D_R(w)$,
\[
D_\Sigma(ws) = \{\alpha_s\}\cup \{s(\beta) \mid \beta\in
D_\Sigma(w)\text{ and }-1<(\beta\mid \alpha_s)<1\}.
\]
The recognising automaton is constructed in the following way. Its
vertex set (which turns out to be always finite) is
\[
\{D_\Sigma(w)\mid w\in W\}.
\]
The labelled edges are given by
\[
D_\Sigma(w) \buildrel s \over \longrightarrow D_\Sigma(ws),
\]
whenever $s\not \in D_R(w)$. Our distinguished starting vertex is
$D_\Sigma(e) = \emptyset$.

\section{Affine and not locally finite groups}\label{se:general}
A {\em well-partially-ordered} (wpo) set is a poset in which every
non-empty subset 
has a minimal element and every antichain is finite. The origin
of the following easy lemma is non-trivial to establish. See e.g.\
Kruskal's survey \cite{kruskal}. We include a proof for convenience
and completeness. 
\begin{lemma} \label{le:finiteproduct}
Suppose $P$ and $Q$ are wpo posets. Then the product poset $P\times Q$
is also wpo. 
\begin{proof}
It is easy to see that the non-empty subsets of $P\times Q$ have
minimal elements. Suppose, in order to
deduce a contradiction, that $A=\{(p_i,q_i)\}_{i\in \N}$ is an
infinite antichain in $P\times Q$. We may assume that $\{p_i\}_{i\in
  \N}$ and $\{q_i\}_{i\in \N}$ are infinite; otherwise we could find
an infinite subset $B\subseteq A$ isomorphic to a subposet of $Q$ or $P$,
respectively, giving a contradiction.

It follows from Ramsey's Theorem that every infinite poset either has
an infinite chain or an infinite antichain (or both). The
antichains in $P$ are finite, so the set $\{p_i\}_{i\in \N}$ contains an
infinite chain. This chain has a smallest element since $P$ is
wpo. Without loss of generality we may therefore assume
$p_0<p_1<p_2<\dots$. Similarly, we may assume that
$\{q_i\}_{i\in \N}$ forms an infinite chain in $Q$ . Since $Q$ is
wpo, we cannot have $q_0>q_1>q_2>\dots$. Therefore, there
exist indices $i$ and $j$ such that $p_i<p_j$ and $q_i\leq q_j$,
contradicting the fact that $A$ is an antichain.
\end{proof}
\end{lemma}
Since non-empty subsets of Coxeter groups always contain minimal
elements under weak order, Lemma
\ref{le:finiteproduct} allows us to restrict attention to irreducible
Coxeter groups --- a group contains an infinite antichain if and only if one
of its irreducible components does.
\begin{theorem} \label{th:affine}
Affine Weyl groups have no infinite antichains.
\begin{proof}
Let $W$ be a finite Weyl group with root system $\Phi$ and associated
affine group $\AW$. Consider the realisation of $\AW$ as a group generated by
affine reflections in $V = \Span \Phi$ (see e.g.\ \cite[Section
  4]{humphreys}). Identifying $V$ with its dual, the reflecting
(affine) hyperplanes are given by 
$H_{\alpha,k} = \{\lambda \in V\mid \langle \lambda,\alpha\rangle = k\}$ for
$\alpha \in \Phi^+$, $k\in \Z$. The complement $V\setminus
\cup_{\alpha,k} H_{\alpha,k}$ is a disjoint union of connected open
    {\em alcoves}. The set of alcoves is in bijection with $\AW$. The
    alcove corresponding to $w\in \AW$ is defined by a (possibly redundant) set of inequalities of
    the form $n_w^\alpha < \langle \lambda,\alpha\rangle < n_w^\alpha + 1$, where
    $n_w^\alpha \in \Z$ for all $\alpha \in \Phi^+$. Corresponding to
    the identity $e\in \AW$, the {\em fundamental 
    alcove} is obtained by putting all $n_e^\alpha = 0$.

Define a partial order $\preceq$ on $\Z$ by letting $i\preceq j$ iff $|i|\leq|j|$ and
either $i=0$ or $\sgn(i) = \sgn(j)$. Thus, $\ldots \succ 2\succ 1 \succ
0 \prec -1 \prec -2 \prec \dots$. Let $\ZZ$ denote this poset. It is
known that the weak 
    order on $\AW$ corresponds to inclusion on the sets of hyperplanes
    that separate the various alcoves from the fundamental one. This
    amounts to saying that, choosing some total ordering of $\Phi^+$,
    the map $\varphi: \AW \to \ZZ^{|\Phi^+|}$ given by $w \mapsto
      (n_w^\alpha)_{\alpha \in \Phi^+}$ is a poset automorphism from
      the weak order on $\AW$ to the image of $\varphi$.

By Lemma \ref{le:finiteproduct}, the antichains in $\ZZ^{|\Phi^+|}$
are finite, and the theorem follows.
\end{proof}
\end{theorem}

A Coxeter group $W$ is called {\em locally finite} if $|W_J|<\infty$
for all $J \subsetneq S$.

\begin{theorem} \label{th:notlocfin}
If $W$ is irreducible and not locally finite, then it has an infinite
antichain. 
\begin{proof}
Suppose $W_J$ is infinite and irreducible for $J\subsetneq S$. Choose
$s\in J$ and $s^\prime\in S\setminus J$ which are neighbours in the
Coxeter diagram of $W$, 
i.e.\ $s$ and $s^\prime$ do not commute. It follows from \cite[Proposition 4.2]{deodhar} that $W_J^{J\setminus \{s\}}$ is infinite. Observe
that $D_R(w)=\{s\}$ for all $w \in W_J^{J\setminus \{s\}}$. By the Word
Property, since $s$ and $s^\prime$ do not commute, every reduced
expression for $ws^\prime$, $w \in W_J^{J\setminus \{s\}}$, contains
exactly one $s^\prime$, and this is necessarily the last letter. This
implies that the infinite set $\{ws^\prime \mid w\in W_J^{J\setminus
  \{s\}}\}$ is an antichain under weak order.
\end{proof}
\end{theorem}

\section{Compact hyperbolic groups}\label{se:comphyp}
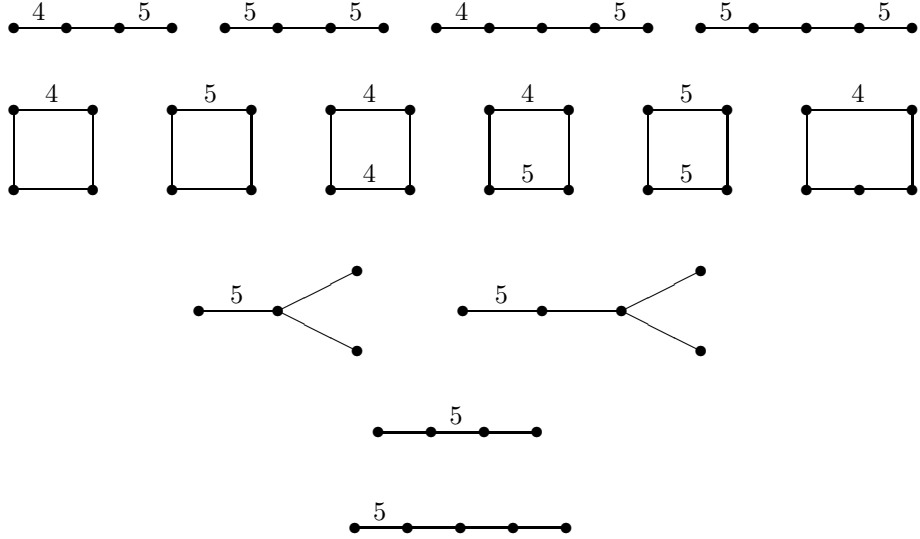
\begin{figure}[htb]
\begin{picture}(360,50)(-10,40) 
\put(0,50){\circle*{4}}
\put(20,50){\circle*{4}}
\put(40,50){\circle*{4}}
\put(60,50){\circle*{4}}
\put(0,50){\line(1,0){20}}
\put(20,50){\line(1,0){20}}
\put(40,50){\line(1,0){20}}
\put(7,53){$4$}
\put(47,53){$5$}

\put(80,50){\circle*{4}}
\put(100,50){\circle*{4}}
\put(120,50){\circle*{4}}
\put(140,50){\circle*{4}}
\put(80,50){\line(1,0){20}}
\put(100,50){\line(1,0){20}}
\put(120,50){\line(1,0){20}}
\put(87,53){$5$}
\put(127,53){$5$}

\put(160,50){\circle*{4}}
\put(180,50){\circle*{4}}
\put(200,50){\circle*{4}}
\put(220,50){\circle*{4}}
\put(240,50){\circle*{4}}
\put(160,50){\line(1,0){20}}
\put(180,50){\line(1,0){20}}
\put(200,50){\line(1,0){20}}
\put(220,50){\line(1,0){20}}
\put(167,53){$4$}
\put(227,53){$5$}

\put(260,50){\circle*{4}}
\put(280,50){\circle*{4}}
\put(300,50){\circle*{4}}
\put(320,50){\circle*{4}}
\put(340,50){\circle*{4}}
\put(260,50){\line(1,0){20}}
\put(280,50){\line(1,0){20}}
\put(300,50){\line(1,0){20}}
\put(320,50){\line(1,0){20}}
\put(267,53){$5$}
\put(327,53){$5$}
\end{picture}
\begin{picture}(360,60)(-10,40) 
\put(0,50){\circle*{4}}
\put(30,50){\circle*{4}}
\put(0,80){\circle*{4}}
\put(30,80){\circle*{4}}
\put(0,50){\line(1,0){30}}
\put(0,80){\line(1,0){30}}
\put(0,50){\line(0,1){30}}
\put(30,50){\line(0,1){30}}
\put(12,83){$4$}

\put(60,50){\circle*{4}}
\put(90,50){\circle*{4}}
\put(60,80){\circle*{4}}
\put(90,80){\circle*{4}}
\put(60,50){\line(1,0){30}}
\put(60,80){\line(1,0){30}}
\put(60,50){\line(0,1){30}}
\put(90,50){\line(0,1){30}}
\put(72,83){$5$}

\put(120,50){\circle*{4}}
\put(150,50){\circle*{4}}
\put(120,80){\circle*{4}}
\put(150,80){\circle*{4}}
\put(120,50){\line(1,0){30}}
\put(120,80){\line(1,0){30}}
\put(120,50){\line(0,1){30}}
\put(150,50){\line(0,1){30}}
\put(132,83){$4$}
\put(132,53){$4$}

\put(180,50){\circle*{4}}
\put(210,50){\circle*{4}}
\put(180,80){\circle*{4}}
\put(210,80){\circle*{4}}
\put(180,50){\line(1,0){30}}
\put(180,80){\line(1,0){30}}
\put(180,50){\line(0,1){30}}
\put(210,50){\line(0,1){30}}
\put(192,83){$4$}
\put(192,53){$5$}

\put(240,50){\circle*{4}}
\put(270,50){\circle*{4}}
\put(240,80){\circle*{4}}
\put(270,80){\circle*{4}}
\put(240,50){\line(1,0){30}}
\put(240,80){\line(1,0){30}}
\put(240,50){\line(0,1){30}}
\put(270,50){\line(0,1){30}}
\put(252,83){$5$}
\put(252,53){$5$}

\put(300,50){\circle*{4}}
\put(320,50){\circle*{4}}
\put(340,50){\circle*{4}}
\put(300,80){\circle*{4}}
\put(340,80){\circle*{4}}
\put(300,50){\line(0,1){30}}
\put(340,50){\line(0,1){30}}
\put(300,50){\line(1,0){20}}
\put(320,50){\line(1,0){20}}
\put(300,80){\line(1,0){40}}
\put(317,83){$4$}
\end{picture}

\begin{picture}(200,50)(0,35) 
\put(0,50){\circle*{4}}
\put(30,50){\circle*{4}}
\put(60,65){\circle*{4}}
\put(60,35){\circle*{4}}
\put(0,50){\line(1,0){30}}
\put(30,50){\line(2,1){30}}
\put(30,50){\line(2,-1){30}}
\put(12,53){$5$}

\put(100,50){\circle*{4}}
\put(130,50){\circle*{4}}
\put(160,50){\circle*{4}}
\put(190,65){\circle*{4}}
\put(190,35){\circle*{4}}
\put(100,50){\line(1,0){30}}
\put(130,50){\line(1,0){30}}
\put(160,50){\line(2,1){30}}
\put(160,50){\line(2,-1){30}}
\put(112,53){$5$}
\end{picture}

\begin{picture}(200,40)(12,40) 
\put(80,50){\circle*{4}}
\put(100,50){\circle*{4}}
\put(120,50){\circle*{4}}
\put(140,50){\circle*{4}}
\put(80,50){\line(1,0){20}}
\put(100,50){\line(1,0){20}}
\put(120,50){\line(1,0){20}}
\put(107,53){$5$}
\end{picture}

\begin{picture}(200,45)(21,30) 
\put(80,50){\circle*{4}}
\put(100,50){\circle*{4}}
\put(120,50){\circle*{4}}
\put(140,50){\circle*{4}}
\put(160,50){\circle*{4}}
\put(80,50){\line(1,0){20}}
\put(100,50){\line(1,0){20}}
\put(120,50){\line(1,0){20}}
\put(140,50){\line(1,0){20}}
\put(87,53){$5$}
\end{picture}
\caption{All irreducible compact hyperbolic Coxeter groups of rank at
  least $4$.}\label{fi:comphyp}
\end{figure}

Lann\'er \cite{lanner} showed that the locally finite Coxeter groups that are
neither finite nor affine are precisely the compact hyperbolic
ones. In rank $3$, every infinite, non-affine group is compact
hyperbolic. The diagrams of the remaining irreducible compact
hyperbolic Coxeter groups are shown in Figure \ref{fi:comphyp}. 

In light of Theorems \ref{th:affine} and \ref{th:notlocfin}, the next
result concludes the proof of Theorem \ref{th:main}.
\begin{theorem}\label{th:comphyp}
Every compact hyperbolic Coxeter group has an infinite antichain.
\end{theorem}
Proving this theorem is the topic of the remainder of the
paper. Proceeding in a case-by-case fashion, the proof is somewhat 
unsatisfactory. In particular, our argument that the group at the bottom
of Figure \ref{fi:comphyp} has an infinite antichain relies on
computer aided calculations and the structure of the automaton
discussed in Section \ref{se:prel}. It would be very interesting 
to have a type-independent proof of Theorem \ref{th:comphyp}, perhaps
in terms of general properties of the symmetric bilinear form
$(\cdot \mid \cdot)$; see
\cite[Section 6.8]{humphreys} for details.

The following simple lemma turns out to produce infinite
antichains in most compact hyperbolic groups.
\begin{lemma}\label{le:antichain}
Suppose that $u,w\in W$ fulfil the following requirements:
\begin{itemize}
\item[(i)] $\ell(u)\leq\ell(w)$.
\item[(ii)] $u \not \leq_R w$.
\item[(iii)] $|S(w)| \geq 3$.
\item[(iv)] Every reduced expression for $wu$ is a concatenation of a
  reduced expression for $w$ and a reduced expression for $u$.
\item[(v)] Every reduced expression for $w^2$ is a concatenation of
  two reduced expressions for $w$.
\end{itemize}
Then, $\{w^ku\mid k\in \N\}$ is an infinite antichain in $W$.
\begin{proof}
Suppose $u$ and $w$ satisfy the hypotheses. We claim that every
reduced expression for $w^ku$, $k\in \N$, is a concatenation of $k$ reduced
expressions for $w$ and one for $u$. To see this, take an
expression for $w^ku$ of the form just described. By (iii), it allows
no braid move which involves an entire copy of $w$. Thus, (iv) and (v)
imply that every braid move simply replaces one copy of $w$ (or $u$)
with another. Moreover, nil moves cannot be possible since it
would mean either that $\ell(wu)<\ell(w)+\ell(u)$ (contradicting (iv))
or that $\ell(w^2)<2\ell(w)$ (contradicting (v)). The claim is proved.

Now assume $w^ku <_R w^lu$ for some $k < l$. By (i) and the above
claim, this means that some reduced expression for $w$ has a reduced
expression for $u$ as a prefix, contradicting (ii). We conclude that
$\{w^ku\}$ is indeed an antichain. 
\end{proof}
\end{lemma}
We say that $u$ and $w$ form a {\em good pair} if they satisfy the
hypotheses of Lemma \ref{le:antichain}.

\begin{lemma}\label{le:increase}
Suppose $W^\prime$ is a Coxeter group obtained from $W$ by increasing
some edge labels in the Coxeter diagram. If $W$ has infinite
antichains, then so does $W^\prime$.
\begin{proof}
Take an infinite antichain $\{w_1, w_2, \dots\}\subset W$. Pick reduced expressions for the $w_i$. These expressions are reduced
in $W^\prime$, too. This is because any sequence of braid moves applicable in
the context of $W^\prime$ is also applicable in $W$; otherwise the
expression would not be reduced in $W$. Thus, the sequence never leads to
a nil move. The corresponding elements therefore form an antichain in
$W^\prime$, too.
\end{proof}
\end{lemma}

\begin{proof}[Proof of Theorem \ref{th:comphyp}] The proof is divided
  into six different cases that combine to exhaust all irreducible compact
  hyperbolic groups (after allowing edge labels to increase,
  using Lemma \ref{le:increase}). The rank three groups are covered by Cases
  I--III. The remaining groups are those in Figure
  \ref{fi:comphyp}. The first row is covered by Case II, the second by
  Case I and the third by Case IV. Finally, the singleton fourth and
  fifth rows are covered by Cases V and VI, respectively. In each
  case except the last one, we apply Lemma \ref{le:antichain} by producing
  a good pair of elements in the corresponding group.

\noindent {\bf Case I.} 
\begin{picture}(200,30)(-40,0)
\put(20,0){\circle*{4}}
\put(60,0){\circle*{4}}
\put(0,20){\circle*{4}}
\put(20,20){\circle*{4}}
\put(60,20){\circle*{4}}
\put(80,20){\circle*{4}}
\put(20,0){\line(-1,1){20}}
\put(60,0){\line(1,1){20}}
\put(20,0){\line(1,0){40}}
\put(0,20){\line(1,0){20}}
\put(60,20){\line(1,0){20}}
\put(20,20){\line(1,0){10}}
\put(60,20){\line(-1,0){10}}
\put(34,19){$\dots$}
\put(37,2){$4$}
\end{picture}

Suppose the Coxeter generators form a cycle $s_1, \dots, s_n, s_1$ in the
Coxeter diagram, and assume $m(s_1,s_n) = 4$. Then, $s_n$ and $s_2\dots
s_ns_1$ form a good pair.

\noindent {\bf Case II.}
\begin{picture}(200,30)(-16,-3)
\put(20,0){\circle*{4}}
\put(40,0){\circle*{4}}
\put(60,0){\circle*{4}}
\put(100,0){\circle*{4}}
\put(120,0){\circle*{4}}
\put(140,0){\circle*{4}}
\put(20,0){\line(1,0){20}}
\put(40,0){\line(1,0){20}}
\put(60,0){\line(1,0){10}}
\put(100,0){\line(-1,0){10}}
\put(100,0){\line(1,0){20}}
\put(120,0){\line(1,0){20}}
\put(74,-1){$\dots$}
\put(27,3){$5$}
\put(127,3){$4$}
\end{picture}

Assume $s_1, \dots, s_n$ is a path in the Coxeter diagram, with
$m(s_1, s_2) = 5$ and $m(s_{n-1}, s_n) = 4$. Now, $s_1s_2s_1$ and $s_1\dots
s_ns_{n-1}\dots s_2$ yield a good pair.

\noindent {\bf Case III.}
\begin{picture}(200,30)(-11,-3)
\put(20,0){\circle*{4}}
\put(40,0){\circle*{4}}
\put(60,0){\circle*{4}}
\put(20,0){\line(1,0){20}}
\put(40,0){\line(1,0){20}}
\put(27,3){$7$}
\end{picture}

Suppose we have three generators $s,t,u\in S$ with $m(s,t) = 7$ and
$m(t,u)\geq 3$. A good pair is given by $st$ and $sutst$.

\noindent{\bf Case IV.}
\begin{picture}(200,30)(-32,-10) 
\put(0,0){\circle*{4}}
\put(30,0){\circle*{4}}
\put(60,0){\circle*{4}}
\put(100,0){\circle*{4}}
\put(130,10){\circle*{4}}
\put(130,-10){\circle*{4}}
\put(0,0){\line(1,0){30}}
\put(30,0){\line(1,0){30}}
\put(60,0){\line(1,0){10}}
\put(100,0){\line(-1,0){10}}
\put(100,0){\line(3,-1){30}}
\put(100,0){\line(3,1){30}}
\put(74,-1){$\dots$}
\put(11,2){$5$}
\end{picture}

Let $s_1, \dots, s_{n-1}$ be a path in the diagram with
$m(s_1,s_2)=5$, and add $s_n$ 
with the relation $m(s_{n-2},s_n)=3$. We get a ``fork-shaped'' diagram
like the one depicted. In this case, we may
recycle the solution from Case II above, regarding the product of the
``ends of the fork'' as a single generator. We get a good pair
consisting of the elements $s_1s_2s_1$ and $s_1\dots
s_ns_{n-2}s_{n-3}\dots s_2$. 

\noindent {\bf Case V.} 
\begin{picture}(200,30)(-36,-3) 
\put(0,0){\circle*{4}}
\put(20,0){\circle*{4}}
\put(40,0){\circle*{4}}
\put(60,0){\circle*{4}}
\put(0,0){\line(1,0){20}}
\put(20,0){\line(1,0){20}}
\put(40,0){\line(1,0){20}}
\put(27,3){$5$}
\end{picture}

Suppose $S=\{s,t,u,v\}$ with $m(s,t)=m(u,v)=3$ and
$m(t,u)=5$. Consider the element $\omega = utvsut = utvust = uvtust =
utsvut = uvtsut$, and note that these are all the reduced expressions for
$\omega$. One easily checks that $\omega$ satisfies condition (v) of
Lemma \ref{le:antichain} by inspecting the 25 concatenations of two
such expressions, observing that none of them
admits a braid move involving both copies of $\omega$. Similarly, we
make sure that $\omega$ together with $\nu = uvtut = utvut$ obey
condition (iv). Conditions (i)--(iii) are immediate, implying that
$(\nu, \omega)$ is a good pair.

\noindent {\bf Case VI.}
\begin{picture}(200,30)(-32,-3) 
\put(0,0){\circle*{4}}
\put(20,0){\circle*{4}}
\put(40,0){\circle*{4}}
\put(60,0){\circle*{4}}
\put(80,0){\circle*{4}}
\put(0,0){\line(1,0){20}}
\put(20,0){\line(1,0){20}}
\put(40,0){\line(1,0){20}}
\put(60,0){\line(1,0){20}}
\put(7,3){$5$}
\end{picture}

Finally, we assume $S = \{s,t,u,v,w\}$, $m(s,t)=5$ and
$m(t,u)=m(u,v)=m(v,w)=3$. Let $\alpha = stuvwstuv$. We claim that
\[
\{\alpha^kw \mid k\in \N\text{ and }k\equiv 0 \pmod{6}\}
\]
 is an infinite antichain. Aided by the computer, we
have established the following facts:
\begin{enumerate}
\item $D_\Sigma(w\alpha^6) = D_\Sigma(w\alpha^7)$. \label{fact1}
\item $\ell(w\alpha^7w) = 2 + 7\ell(\alpha) = 65$. \label{fact2}
\end{enumerate} 

Now consider the recognising automaton for the language of reduced
expressions described in Section \ref{se:prel}. Combining (\ref{fact1})
and (\ref{fact2}), we see that there is a cycle of length $\ell(\alpha)=9$
corresponding to $\alpha$ which begins and ends in
$D_\Sigma(w\alpha^6)$. Repeatedly traversing this cycle produces
reduced expressions. Furthermore, fact (\ref{fact2}) shows that after
walking this cycle a number of times, we may use an edge labelled $w$. Thus,
$\ell(w\alpha^kw)=2+k\ell(\alpha)$ for all $k\geq 6$. This proves
$w\not <_R\alpha^kw$ for such $k$. Since $\alpha^kw
<_R\alpha^lw \Leftrightarrow w<_R\alpha^{l-k}w$, for $k<l$, we
conclude that $\{\alpha^kw\mid k\in 6\N\}$ is indeed an infinite antichain.

\end{proof}

\end{document}